\documentclass[11pt]{amsart}
\usepackage{amssymb,mathrsfs}
\title{A Note on the Artin Conjecture}

\author{Jae-Hyun Yang}
\address{Department of Mathematics, Inha University,
Incheon 402-751, Korea}
\email{jhyang@inha.ac.kr }

\begin{document}

\newtheorem{theorem}{Theorem}[section]
\newtheorem{lemma}{Lemma}[section]
\newtheorem{proposition}{Proposition}[section]
\newtheorem{remark}{Remark}[section]
\newtheorem{definition}{Definition}[section]
\newtheorem{Cojecture}{Conjecture}[section]

\renewcommand{\theequation}{\thesection.\arabic{equation}}
\renewcommand{\thetheorem}{\thesection.\arabic{theorem}}
\renewcommand{\thelemma}{\thesection.\arabic{lemma}}
\newcommand{\BR}{\mathbb R}
\newcommand{\BQ}{\mathbb Q}
\newcommand{\bn}{\bf n}
\def\charf {\mbox{{\text 1}\kern-.24em {\text l}}}
\newcommand{\BC}{\mathbb C}
\newcommand{\BZ}{\mathbb Z}


%\maketitle

\begin{abstract}
{In this paper, we survey some recent results on the Artin
conjecture and discuss some aspects for the Artin conjecture.}
\end{abstract}
\maketitle

%\tableofcontents
\newcommand\tr{\triangleright}
\newcommand\al{\alpha}
\newcommand\be{\beta}
\newcommand\g{\gamma}
\newcommand\gh{\Cal G^J}
\newcommand\G{\Gamma}
\newcommand\de{\delta}
\newcommand\e{\epsilon}
\newcommand\z{\zeta}
\newcommand\vth{\vartheta}
\newcommand\vp{\varphi}
\newcommand\om{\omega}
\newcommand\p{\pi}
\newcommand\la{\lambda}
\newcommand\lb{\lbrace}
\newcommand\lk{\lbrack}
\newcommand\rb{\rbrace}
\newcommand\rk{\rbrack}
\newcommand\s{\sigma}
\newcommand\w{\wedge}
\newcommand\fgj{{\frak g}^J}
\newcommand\lrt{\longrightarrow}
\newcommand\lmt{\longmapsto}
\newcommand\lmk{(\lambda,\mu,\kappa)}
\newcommand\Om{\Omega}
\newcommand\ka{\kappa}
\newcommand\ba{\backslash}
\newcommand\ph{\phi}
\newcommand\M{{\Cal M}}
\newcommand\bA{\bold A}
\newcommand\bH{\bold H}

\newcommand\Hom{\text{Hom}}

\newcommand\pa{\partial}

\newcommand\pis{\pi i \sigma}
\newcommand\sd{\,\,{\vartriangleright}\kern -1.0ex{<}\,}
\newcommand\wt{\widetilde}
\newcommand\fg{\frak g}
\newcommand\fk{\frak k}
\newcommand\fp{\frak p}
\newcommand\fs{\frak s}
\newcommand\fh{\frak h}
\newcommand\Cal{\mathcal}

\newcommand\fn{{\frak n}}
\newcommand\fa{{\frak a}}
\newcommand\fm{{\frak m}}
\newcommand\fq{{\frak q}}
\newcommand\CP{{\mathcal P}_g}
\newcommand\BD{\mathbb D}
\newcommand\BH{\mathbb H}
\newcommand\CM{{\mathcal M}}

\newcommand\OF{\overline F}
\newcommand\GF{\textrm{Gal}({\overline F}/F)}
\newcommand\OFv{\overline{F_v}}
\newcommand\GQ{ \textrm{Gal}({\overline \BQ}/\BQ)}

%%%%%%%%%%%%%%%%%%%%%%%%%%%%%%%%%%%%%%%%%%%%%%%%%%%%%%%%%%%%%%%%%%%
%
%                             Sec 1  Introduction
%
%%%%%%%%%%%%%%%%%%%%%%%%%%%%%%%%%%%%%%%%%%%%%%%%%%%%%%%%%%%%%%%%%%%

\begin{section}{{\bf Introduction}}
\setcounter{equation}{0} Let $K/\BQ$ be a Galois extension of
$\BQ$ and $\rho: \textrm{Gal}(K/\BQ)\lrt GL(n,\BC)$ a nontrivial
irreducible representation of its Galois group. E. Artin [1]
associated to this data an $L$-function $L(s,\rho)$, defined for $
\textrm{Re}\,s>1,$ which he conjectured to continue analytically
to an entire function on the whole complex plane $\BC$ satisfying
a functional equation. In 1947, R. Brauer\,[6] showed that the
Artin $L$-function $L(s,\rho)$ has a meromorphic continuation to a
meromorphic function on $\BC$ and satisfies a functional
equation.\\ \indent
%\vskip 0.1cm
Artin established his conjecture for the $ \textit{monomial}$
representations, those induced from one-dimensional representation
of a subgroup. His conjecture has not been solved yet in any
dimension $\geq 2$. More evidence is provided in dimension $2$ by
R. Langlands, J. Tunnell, R. Taylor et al. In the case of two
dimensional icosahedral representations, his conjecture still
remains open. When $\rho$ is an odd icosahedral representation,
infinitely many examples of the Artin conjecture are known by the
work of R. Taylor and others.\\ \indent
%\vskip 0.2cm
This article is organized as follows. In Section 2, we review the
Galois representations of $\GQ$ roughly. In Section 3, we describe
the definition of the Artin $L$-function. In Section 4, we explain
the connection between the Artin conjecture and the Langlands
Functoriality Conjecture. In Section 5, we survey some known
results on the Artin conjecture in the two dimensional case. In
the final section we discuss some aspects for the Artin
conjecture. \vskip 0.1cm \noindent $\textsc{Notations\,:}$
Throughout this paper, $F$ denotes a number field, ${\overline F}$
an algebraic closure of $F$, and $ \textrm{Gal}({\overline F}/F)$
the absolute Galois group of $F$. We regard $\GF$ as a topological
group relative to the Krull topology. We write ${\mathbb A}_F$ and
$I_F$ for the adele ring and the idele group attached to $F$
respectively. For each place $v$ of $F$, we let $F_v$ be the
completion of $F$ relative to $v$. We also fix an algebraic
closure ${\overline{F_v}}$ of $F_v$ for each place $v$. For a
square matrix $A$, $ \textrm{tr}(A)$ denotes the trace of $A$.

\end{section}

\newcommand\BF{\mathbb F}
\newcommand\OQ{\overline{\BQ}}
\newcommand\Qp{\BQ_p}
\newcommand\OQp{\overline{\BQ}_p}
\newcommand\Fp{{\mathbb F}_p}
\newcommand\OFp{{\overline{\mathbb F}}_p}
\newcommand\GQp{ \textrm{Gal}({\overline{\BQ}}_p/\BQ_p) }
\newcommand\GQl{ \textrm{Gal}({\overline{\BQ}}_l/\BQ_l) }
\newcommand\GFp{ \textrm{Gal}({\overline{\mathbb F}}_p/\BF_p) }
\newcommand\Zp{\BZ_p}
\newcommand\COp{{\mathcal O}_{\OQp}}
\newcommand\MOp{{\mathfrak m}_{\OQp}}
\newcommand\Frp{ \textrm{Fr}_p}
\newcommand\OQl{\overline{\BQ}_l}
\newcommand\Wp{W_{\Qp}}

\vskip 0.7cm
%%%%%%%%%%%%%%%%%%%%%%%%%%%%%%%%%%%%%%%%%%%%%%%%%%%%%%%%%%%%%%%%%%%
%
%   Sec 2  Galois Representations
%
%%%%%%%%%%%%%%%%%%%%%%%%%%%%%%%%%%%%%%%%%%%%%%%%%%%%%%%%%%%%%%%%%%%
\begin{section}{{\bf Galois Representations}}
\setcounter{equation}{0} R. Taylor published a good survey paper
[23] about Galois representations. The content of this section is
a brief description of Section 1 in [23].\par Let $\BQ$ be the
field of rational numbers and $\OQ$ denote the algebraic closure
of $\BQ$. We let $\GQ$ be the absolute Galois group of $\BQ$. We
see that $\GQ$ is a profinite topological group, a basis of open
neighborhoods of the identity being given by the subgroups $
\textrm{Gal}(\OQ/K)$ as $K$ runs over subextensions of $\OQ/\BQ$
which is finite over $\BQ$. Let $\Qp$ be the field of $p$-adic
numbers, which is a totally disconnected locally compact
topological field. $\OQp/\Qp$ is an infinite extension of $\Qp$
and $\OQp$ is not complete. We shall denote its completion by
$\BC_p$. Let $\Zp$\,(resp.\,$\COp$) be the ring of integers in
$\Qp$\,(resp.\,$\OQp$). These are local rings with maximal ideals
$p\Zp$ and $\MOp$ respectively. Then it is easy to see that the
field $\OFp:=\COp/\MOp$ is an algebraic closure of the field
$\Fp:=\Zp/p\Zp.$ Thus we obtain a continuous map
\begin{equation*}
\GQp\lrt \GFp
\end{equation*}

\noindent which is surjective. Its kernel is called the $
\textit{inertia subgroup}$ of $\GQp$, and is denoted by $I_{\Qp}$.
The Galois group $\GFp$ is procyclic and has a canonical generator
$\Frp$ called the Frobenius element defined by
\begin{equation*}
\Frp(x)=x^p,\quad x\in \OFp.
\end{equation*}

I want to describe $\GQ$ via its representations. We have two
natural representations of $\GQ$, which are
\begin{equation*}
\GQ\lrt GL(n,\BC),\quad \textrm{the Artin representations}
\end{equation*}
and
\begin{equation*}
\GQ\lrt GL(n,\OQl),\quad \textrm{the}\ l\textrm{-adic
representations}.
\end{equation*}

\noindent Here $GL(n,\OQl)$ is a group with $l$-adic topology.
These representations are continuous. \vskip 0.2cm  The $l$-adic
representations are closely related to an arithmetic geometry.
%\begin{enumerate}
\begin{itemize}
\item A choice of embeddings $\OQ\hookrightarrow \BC$ and
$\OQ\hookrightarrow\OQl$ establishes a bijection between
isomorphism classes of Artin representations and isomorphism
classes of $l$-adic representations with open kernel.
\item There is a unique character
\begin{equation*}
\chi_l:\GQ\lrt \BZ_l^{\times}\subset \OQ_l^{\times}
\end{equation*}
such that
\begin{equation*}
\s \zeta=\zeta^{\chi_l(\s)}
\end{equation*}
for all $l$-power roots of unity $\zeta$. This is called the
$l$-$\textit{adic cyclotomic character.}$
\item If $X/\BQ$ is a smooth projective variety, then the natural
action of $\GQ$ on the cohomology
\begin{equation*}
H^i(X(\BC),\OQl)\cong H_{et}^i(X\times_{\BQ}\OQ,\OQl)
\end{equation*}
is an $l$-adic representation.
\end{itemize}
%\end{enumerate}
\vskip 0.1cm\ \ We now discuss $l$-adic representations of $\GQp$.
Let $W_{\Qp}$ be the subgroup of $\GQp$ consisting of elements
$\s\in \GQp$ such that $\s$ maps to $\Frp^{\BZ}\subset \GFp.$ We
endow $W_{\Qp}$ with a topology by decreeing that $I_{\Qp}$ with
its usual topology should be an open subgroup of $W_{\Qp}$. We
first consider the case $l\neq p.$ We define a
$WD$-$\textit{representation}$ of $\Wp$ over a field $E$ to be a
pair
\begin{equation*}
r:\Wp\lrt GL(V),\quad \textrm{a continuous representation of $\Wp$
with open kernel }
\end{equation*}
and
\begin{equation*}
N\in \textrm{End}(V),\quad \textrm{a nilpotent endomorphism of
$V$}
\end{equation*}
\noindent such that
\begin{equation*}
r(\phi)\,N\,r(\phi^{-1})=p^{-1}N
\end{equation*}

\noindent for every lift $\phi\in \Wp$ of $\Frp$, where $V$ is a
finite dimensional $E$-vector space. A WD-representation $(r,N)$
is said to be $ \textit{unramified}$ if $N=0$ and
$r(I_{\Qp})=\left\{ 1\right\}.$ In the case $E=\OQl$, we call a
WD-representation $(r,N)\ l$-$\textit{integral}$ if all
eigenvalues of $r(\phi)$ has the absolute value $1$. If $l\neq p$,
then there is an equvalence of categories between $l$-integral
WD-representations of $\Wp$ over $\OQl$ and $l$-adic
representations of $\GQp$. We will write $WD_p(R)$ for the
WD-representation associated to an $l$-adic representation $R$ of
$\GQp$. An $l$-adic representation $R$ is said to be $
\textit{unramified}$ if $WD_p(R)$ is unramified. The case $l=p$ is
much more complicated because there are many more $p$-adic
representations of $\GQp$. These have been extensively studied by
J.-M. Fontaine et al. They single out certain special $p$-adic
representations which are called $ \textit{de Rham}$. Indeed most
$p$-adic representations are not de Rham. To any de Rham
representation $R$ of $\GQp$ on a $\OQp$-vector space $V$ they
associate the following pair\,:
\begin{itemize}
\item A WD-representation $\textrm{WD}_p(R)$ of $\Wp$ over $\OQp$.
\item A multiset HT($R$) of $\dim V$ integers, called the
Hodge-Tate numbers of $R$. The multiplicity of $i$ in HT($R$) is
\begin{equation*}
\dim_{\OQp}(V\otimes_{\Qp}\BC_p(i))^{\GQp},
\end{equation*}
where $\BC_p(i)$ denotes $\BC_p$ with $\GQp$-action and $\GQp$
acts on $\BC_p$ via $\chi_p(\s)^i$ times its usual Galois action
on $\BC_p.$
\end{itemize}
\vskip 0.1cm We refer to [10,\,11,\,12] and [2] for more details
on de Rham representations and their related materials. \vskip
0.2cm We now discuss a so-called $\textit{geometric}\ l$-adic
representations. Fontaine and Mazur\,[13] proposed the following
conjecture. \vskip 0.1cm \noindent {\bf Conjecture A}
(Fontaine-Mazur) $ \textit{Suppose that}$
\begin{equation*} R:\GQ\lrt GL(V)
\end{equation*}

\noindent $ \textit{is an irreducible $l$-adic representation
which is unramified at all but finitely many primes and}$\\ $
\textit{with $R|_{\GQl}$ de Rham.\ Then there is a smooth
projective variety $X/\BQ$ and integers $i\geq 0$ }$ \\
$ \textit{and $j$ such that $V$ is a subquotient of
$H^i(X(\BC),\OQl(j))$. In particular $R$ is pure of some wight}$\\
$ \textit{$w\in \BZ$.}$ \vskip 0.1cm Tate formulated the following
conjecture. \vskip 0.1cm\noindent {\bf Conjecture B} (Tate). $
\textit{Suppose that $X/\BQ$ is a smooth projective variety. Then
there is a}$\\ $ \textit{decomposition}$
\begin{equation*}
H^i(X(\BC),\OQ)=\oplus_j M_j
\end{equation*}

\noindent $ \textit{with the following properties\,:}$\\
%$ \textit{
%         \begin{itemize}
%            \begin{enumerate}
%               \item For each prime $l$ and for each embeddings $\iota:\OQ\hookrightarrow \OQl,\ M_j\otimes_{\OQ,l}\OQl$
%               is an irreducible subrepresentation of $H^i(X(\BC),\OQl).$
%               \item For all indices $j$ and for all primes $p$ there is a WD-representation $WD_p(M_j)$ of $\Wp$ over $\OQl$
%                     such that \begin{equation*} WD_p(M_j)\otimes_{\OQl,\iota}\OQl\cong WD_p(M_j\otimes_{\OQ,\iota}\OQl)
%                     \end{equation*} for all primes $l$ and all embeddings $\iota:\OQ\lrt \OQl.$
%               \item There is a multiset of integers $ \textrm{HT}(M_j)$ such that
%                  \begin{itemize}
%                   \begin{enumerate}
%                      \item for all primes $l$ and all embeddings  $\iota:\OQ\hookrightarrow \OQl$
%                             \begin{equation*}    \textrm{HT}(M_j\otimes_{\OQ,\iota}\OQl)= \textrm{HT}(M_j)
%                             \end{equation*}
%                      \item and for all $\iota:\OQ\hookrightarrow \BC$
%                             \begin{equation*}    \dim_{\BC}(M_j\otimes_{\OQ,\iota}\BC)\cap H^{a,i-a}(X(\BC),\BC)
%                             \end{equation*}
%                             is the multiplicity of $a$ in $
%                             \textrm{HT}(M_j)$.
%                    \end{enumerate}
%                    \end{itemize}
%               \end{enumerate}
%             \end{itemize}
%             }$
\vskip 0.1cm\ \ \ $ \textit{1.\ For each prime $l$ and for each
embeddings $\iota:\OQ\hookrightarrow \OQl,\
M_j\otimes_{\OQ,l}\OQl$ is an irreducible}$\par \ \ \ \ \ \
$\textit{subrepresentation of $H^i(X(\BC),\OQl).$}$ \par \ \ \
$\textit{2. For all indices $j$ and for all primes $p$ there is a
WD-representation $WD_p(M_j)$ of}$ \par\ \ \ \ \ \ \ $
\textit{$\Wp$ over $\OQl$ such that}$
\begin{equation*}
WD_p(M_j)\otimes_{\OQ,\iota}\OQl\cong
WD_p(M_j\otimes_{\OQ,\iota}\OQl)
\end{equation*}

\indent\ \ \ \ \ \ \ $ \textit{for all primes $l$ and all
embeddings $\iota:\OQ\lrt \OQl.$}$
\par \ \ \
$\textit{3. There is a multiset of integers $ \textrm{HT}(M_j)$
such that}$ \par \ \ \ \ \ \ \ \ $\textit{ (a) for all primes $l$
and all embeddings  $\iota:\OQ\hookrightarrow \OQl$}$
\begin{equation*}    \textrm{HT}(M_j\otimes_{\OQ,\iota}\OQl)= \textrm{HT}(M_j)
                             \end{equation*}

\par \ \ \ \ \ \ \ \ $\textit{ (b) and for all $\iota:\OQ\hookrightarrow
\BC$}$
\begin{equation*}    \dim_{\BC}(M_j\otimes_{\OQ,\iota}\BC)\cap H^{a,i-a}(X(\BC),\BC)
\end{equation*}

\par \ \ \ \ \ \ \ \ \ \ \ \ \ $\textit{is the multiplicity of $a$ in
$\textrm{HT}(M_j)$.}$

\vskip 0.1cm If one believes conjecture A and B, then geometric
$l$-adic representations should come in compatible families as $l$
varies. There are many ways to make precise the notion of such a
compatible family. See [23] for one of such family.

\end{section}

\newcommand\GFv{\textrm{Gal}(\OFv/F_v)}
\newcommand\Okv{\overline{k_v}}
\newcommand\Gkv{\textrm{Gal}(\Okv/k_v)}
\newcommand\Frv{\textrm{Fr}_v}

\vskip 0.7cm
%%%%%%%%%%%%%%%%%%%%%%%%%%%%%%%%%%%%%%%%%%%%%%%%%%%%%%%%%%%%%%%%%%%
%
%   Sec 3  Artin $L$-Functions
%
%%%%%%%%%%%%%%%%%%%%%%%%%%%%%%%%%%%%%%%%%%%%%%%%%%%%%%%%%%%%%%%%%%%
\begin{section}{{\bf Artin $L$-Functions}}
\setcounter{equation}{0}

Let $F$ be a number field. We let
\begin{equation*}
\sigma: \GF\lrt GL(V)
\end{equation*}

\noindent be a finite dimensional Galois representation over $F$,
where $V$ is a finite dimensional complex vector space. The $
\textit{Artin}\ L$-$ \textit{function}\ L(s,\sigma)$ attached to
the Galois representation $\sigma$ is defined to be an Euler
product
\begin{equation*}
L(s,\sigma)=\prod_v L(s,\s_v),
\end{equation*}

\noindent where $v$ runs over all places of $F$. The local factor
$L(s,\s_v)$ is defined as follows. First we choose an embedding
$i_v:\OF\lrt \OFv$, which gives rise to an embedding of Galois
groups

\begin{equation*}
j_v: \textrm{Gal}(\OFv/F_v)\lrt \GF
\end{equation*}

\noindent via $ \textit{restriction}$.
The composition $\s_v=\s\circ j_v$ is a continuous representation
of $\GFv$. It depends on the choice of an embedding $i_v$, but
different choices of $i_v$ lead to conjugate embeddings ${j_v}$.
So the equivalence class of $\s_v$ is well defined and depends $
\textit{only}$ on $v$.
\vskip 0.31cm In the nonarchimedean case, we let $k_v$ and
${\overline{k_v}}$ denote the residue fields of $F_v$ and $\OFv$
respectively. $\GFv$ acts on $\Okv$ and we have an exact sequence
\begin{equation*}
1\lrt I_v\lrt \GFv \lrt \Gkv\lrt 1,
\end{equation*}

\noindent where $I_v$ is the inertia subgroup. We set $q_v=|k_v|$.
A $ \textit{Frobenius\ element}\ \Frv$ is an element of $\GFv$
whose image in $\Gkv$ is the automorphism
\begin{equation*}
x\mapsto x^{q_v},\quad x\in \Okv.
\end{equation*}

\noindent We note that the action of $\s(\Frv)$ on the subspace
$V^{I_v}$ of inertial invariants in $V$ is independent of the
choice $\Frv$. We define the local factor $L(s,\s_v)$ at $v$ by
\begin{equation*}
L(s,\s_v)=\det \left( 1-q_v^{-s}\s_v(\Frv)|_{V^{I_v}}\right)^{-1}.
\end{equation*}

\vskip 0.2cm
 The Galois representation $\s$ is said to be $
\textit{unramified}$ at $v$ if $\s_v(I_v)={1}$. In this case, the
element $\s_v(\Frv)$ is independent of the choice of $\Frv$. The $
\textit{Frobenius class}$ attached to $v$ is the conjugacy class
$\{ \s_v(\Frv)\}$ of $\s_v(\Frv)$ in $GL(V)$. The Frobenius class
is independent of the choice of an embedding $j_v$ and thus
depends only on $v$. We note that it is a semisimple conjugacy
class, that is, it consists of diagonalizable elements.

\vskip 0.2cm If $v$ is archimedean, then $F_v\cong \BR$ or $\BC$.
In case $F_v\cong \BR$, we get $\GFv\cong
\textrm{Gal}(\BC/\BR)=\left\{ 1,c\right\},$ where $c$ denotes the
complex conjugation. The eigenvalues of $\s_v(c)$ are $\pm 1$. Let
$m_+$\,(resp.\,$m_-$) be the number of $+1$\,(resp.\,$-1$)
eigenvalues of $\s_v(c)$. In this case, we define the local factor
$L(s,\s_v)$ by
\begin{equation*}
L(s,\s_v)=\left( \pi^{-s/2}\Gamma(s/2)\right)^{m_+}\left(
\pi^{-(s+1)/ 2}\Gamma((s+1)/ 2)\right)^{m_-}.
\end{equation*}

If $F_v\cong \BC$, then $\GFv\cong \{1\}$. In this case, we define
\begin{equation*}
L(s,\s_v)=\left( 2(2\pi)^{-s}\Gamma(s)\right)^n,
\end{equation*}

\noindent where $n=\dim_{\BC} V.$ \vskip 0.2cm It is easy to see
that
\begin{equation*}
L(s,\s\oplus \tau)=L(s,\s)L(s,\tau)
\end{equation*}

\noindent for any two Galois representations $\s$ and $\tau$. For
any finite set $S$ of places, we define the $\textit{partial}$
$L$-function $L_S(s,\s)$ by
\begin{equation*}
L_S(s,\s)=\prod_{v\notin S}L(s,\s_v).
\end{equation*}

\noindent We observe that if $\s$ is the trivial representation
and $S$ is the archimedean places, then
\begin{equation*}
L_S(s,\s)=\prod_{v<\infty}\left( 1-q_v^{-s}\right)^{-1}
\end{equation*}

\noindent is nothing but the so-called Dedekind zeta function
$\zeta_F(s)$of $F$.

\end{section}

\vskip 0.7cm
%%%%%%%%%%%%%%%%%%%%%%%%%%%%%%%%%%%%%%%%%%%%%%%%%%%%%%%%%%%%%%%%%%%
%
%   Sec 4 The Artin Conjecture and Functoriality
%
%%%%%%%%%%%%%%%%%%%%%%%%%%%%%%%%%%%%%%%%%%%%%%%%%%%%%%%%%%%%%%%%%%%
\begin{section}{{\bf The Artin Conjecture and Functoriality}}
\setcounter{equation}{0}

Since the eigenvalues of $\s_v(\Frv)$ at each place $v$ are roots
of unity, it is easy to see that the Euler product for $L(s,\s)$
converges absolutely for $ \textrm{Re}\,s >1.$ According to the
works of E. Hecke\,[17], E. Artin\,[1] and R. Brauer\,[6], we
obtain the following theorem.
\begin{theorem} Let $F$ be a number field. Let $\s:\GF\lrt GL(V)$
be a finite dimensional complex Galois representation over $F$.
Then the Artin $L$-function $L(s,\s)$ has a meromorphic
continuation to a meromorphic function on $\BC$. Moreover
$L(s,\s)$ satisfies a functional equation
\begin{equation*}
L(s,\s)=\epsilon (s,\s)L(1-s,\s^*),
\end{equation*}

\noindent where $\epsilon(s,\s)$ is the so-called $
\textit{epsilon factor}$\,(cf.\,[22]) and $\s^*$ denotes the
contragredient representation of $\s$.
\end{theorem}

\noindent $\textbf{Artin Conjecture.}$ $ \textit{If $\s:\GF\lrt
GL(V)$ is a nontrivial irreducible finite dimensi-}$\par\noindent
$\textit{onal complex Galois representation over $F$, then the
Artin $L$-function $L(s,\s)$ can be analy-}$\par\noindent $
\textit{tically continued to an entire function on $\BC$}$.

\vskip 0.3cm Let $F$ be a number field. For a cuspidal
representation $\pi$ of $GL(n,{\mathbb A}_F)$, we can define the $
\textit{automorphic}$ $L$-$ \textit{function} \ L(s,\pi)$ of $\pi$
given by
\begin{equation*}
L(s,\pi)=\prod_v L(s,\pi_v),
\end{equation*}

\noindent where $v$ runs over all places of $F$. The precise
definition of $L(s,\pi_v)$ can be found in [5],\,[14],\,[15] and
[19]. Jacquet and Langlands [18] proved that if $n=2,\ L(s,\pi)$
can be analytically continued to an entire function on the whole
complex plane $\BC$. Godement and Jacquet [16] proved that for any
positive integer $n$, the $L$-function $L(s,\pi)$ can be
analytically continued to an entire function on the whole complex
plane $\BC$.

\vskip 0.3cm R. Langlands proposed the following conjecture in
order to attack the Artin conjecture. \vskip 0.2cm \noindent $
\textbf{Langlands Functoriality Conjecture.}$ $ \textit{ Let F be
a number field. Let}$
\begin{equation*}
 \s:\GF\lrt GL(n,\BC)
\end{equation*}

\noindent $\textit{be a nontrivial irreducible finite dimensional
complex Galois representation over F. Then}$ \\ $\textit{there
exists a cuspidal representation}$ $\pi(\s)$ of $GL(n,{\mathbb
A}_F)$ $ \textit{such that}$
\begin{equation*}
L(s,\s)=L(s,\pi(\s)).
\end{equation*}

\vskip 0.3cm We observe that if Langlands Functoriality Conjecture
is true, by the work of Jacquet, Langlands and Godement, the Artin
conjecture is true. The Langlands Functoriality Conjecture gave
rise to some solutions of the Artin conjecture for an irreducible
two-dimensional Galois representation. \vskip 0.3cm I want to
introduce the recent work of Andrew Booker. Let
\begin{equation*}
 \rho: \textrm{Gal}({\overline {\mathbb Q}}/{\mathbb Q})   \lrt GL(2,\BC)
\end{equation*}

\noindent be an irreducible two-dimensional complex Galois
representation over ${\mathbb Q}$. Booker\,[3] proved the
following result.

\begin{theorem} If $L(s,\rho)$ is not automorphic, then it has
infinitely many poles. In particular, the Artin conjecture for
$\rho$ implies the Langlands Functoriality Conjecture for $\rho$.
\end{theorem}

\end{section}

\vskip 0.7cm
%%%%%%%%%%%%%%%%%%%%%%%%%%%%%%%%%%%%%%%%%%%%%%%%%%%%%%%%%%%%%%%%%%%
%
%   Sec 5  Special Cases  of the Artin Conjecture
%
%%%%%%%%%%%%%%%%%%%%%%%%%%%%%%%%%%%%%%%%%%%%%%%%%%%%%%%%%%%%%%%%%%%
\begin{section}{{\bf Special Cases  of the Artin Conjecture}}
\setcounter{equation}{0} Let $F$ be a number field. Let
\begin{equation*}
 \rho: \GF\lrt GL(2,\BC)
\end{equation*}

\noindent be an irreducible two-dimensional complex Galois
representation over $F$. The adjoint representation of the group
$GL(2,\BC)$ on the Lie algebra ${\mathfrak g}{\mathfrak l}(2,\BC)$
induces the adjoint action of $GL(2,\BC)$ on the three dimensional
Lie algebra ${\mathfrak s}{\mathfrak l}(2,\BC)$ of $2\times 2$
complex matrices of trace zero. We denote this representation by
\begin{equation*}
 \textrm{Ad}: GL(2,\BC)\lrt GL(3,\BC).
\end{equation*}

The symmetric bilinear form $ \textrm{tr}(AB)$ is invariant under
the adjoint action of $GL(2,\BC)$, and the image of $ \textrm{Ad}$
is isomorphic to the complex orthogonal group $SO(3,\BC)$ defined
by this bilinear form $ \textrm{tr}(AB)$. Irreducible
two-dimensional representations are classified according to the
image of $ \textrm{Ad}\circ \rho$ in $SO(3,\BC)$. It is known that
a finite subgroup of $SO(3,\BC)$ is either $ \textit{cyclic,
dihedral}$ or isomorphic to one of the symmetry groups of the
Platonic solids\,: \vskip 0.2cm (1) $ \textit{tetrahedral
group}\cong A_4$\,;\\ \indent (2) $ \textit{octahedral group}\cong
S_4$\,;\\ \indent (3) $ \textit{icosahedral group}\cong A_5$.

\vskip 0.3cm We shall say that $\rho$ is of $ \textit{cyclic,
dihedral, tetrahedral, octahedral, icosahedra type}$ if the image
of $\textrm{Ad}\circ \rho$ in $SO(3,\BC)$ is of the corresponding
type. The Artin conjecture was solved by E. Artin for the cyclic
and dihedral type, by R. Langlands [20] for the tetrahedral type,
and was solved completely by J. Tunnell [25,\,26] for the
octahedral type. Indeed we can show that $\pi(\rho)$ exists if
$\rho$ is of cyclic, dihedral, tetrahedral or octahedral type. We
refer to [21] for sketchy proofs in such these types. The Artin
conjecture for the $ \textit{icosahedral type}$ has not been
solved yet, although it has been verified in few very special
cases\,[7,\,8,\,9,\,24]. \vskip 0.2cm Recently Taylor et al gave
some evidences to the Artin conjecture for the icosahedral type.
K. Buzzard, M. Dickinson, N. Sheherd-Baron and R. Taylor\,[8]
proved that the Artin conjecture is true for certain special odd
icosahedral representations of $\textrm{Gal}({\overline {\mathbb
Q}}/{\mathbb Q})$ by showing that they are modular. I describe
this content explicitly.

\begin{theorem} Suppose that $\rho:\textrm{Gal}({\overline {\mathbb Q}}/{\mathbb
Q})\lrt GL(2,\BC)$ is a continuous representation and that $\rho$
is odd, i.e., the determinant of $\rho(c)$ is -1, where $c$ is the
complex conjugation. If $\rho$ is of icosahedral type, we assume
that \vskip 0.2cm $\bullet$ the projectivised representation
proj$(\rho): \textrm{Gal}({\overline {\mathbb Q}}/{\mathbb Q})\lrt
PGL(2,\BC)$ is unramified at $2$ \\ \indent \ \ and that image of
a Frobenius element at $2$ under proj($\rho)$ has order $3$,
\vskip 0.2cm $\bullet$ and proj$(\rho)$ is unramified at $5$.
\vskip 0.2cm \noindent Then there is a new form of weight one such
that for all primes $p$ the $p$-th Fourier coefficient of $f$
equals the trace of Frobenius at $p$ on the inertia at $p$
covariants of $\rho$. In particular the Artin $L$-function for
$\rho$ is the Mellin transform of a newform of weight one and is
an entire function.
\end{theorem}

Moreover R. Taylor\,[24] proved the following.

\begin{theorem} Let $\rho:\textrm{Gal}({\overline {\mathbb Q}}/{\mathbb
Q})\lrt GL(2,\BC)$ is a continuous representation and that $\rho$
is odd, i.e., the determinant of $\rho(c)$ is -1, where $c$ is the
complex conjugation. If $\rho$ is of icosahedral type, we assume
that the projective image of the inertia group at $3$ has odd
order and the projective image of the decomposition group at $5$
is unramified at $2$. Then $\rho$ is modular and its Artin
$L$-function $L(s,\rho)$ is entire.
\end{theorem}

\end{section}

\vskip 0.7cm
%%%%%%%%%%%%%%%%%%%%%%%%%%%%%%%%%%%%%%%%%%%%%%%%%%%%%%%%%%%%%%%%%%%
%
%   Sec 6  Final Remarks
%
%%%%%%%%%%%%%%%%%%%%%%%%%%%%%%%%%%%%%%%%%%%%%%%%%%%%%%%%%%%%%%%%%%%
\begin{section}{{\bf Final Remarks}}
\setcounter{equation}{0} As mentioned before, we still have no
idea of verifying the Artin conjecture for $n$-dimensional Galois
representations with $n\geq 3.$ In the case of two dimensional
icosahedral representations, the Artin conjecture still remains
open. If the Artin conjecture is true for a certain Galois
representation $\rho$, it might be interesting to find methods for
locating zeros of the Artin $L$-function $L(s,\rho)$. In [4], A.
Booker discusses two methods for locating zeros of $L(s,\rho)$. He
also presents a group-theoretic criterion under which one may
verify the Artin conjecture for some non-monomial Galois
representations, up to finite height in the complex plane.

\end{section}

\end{document}